\documentclass[twoside,leqno,11pt]{article}
\input epsfx

\input prepictex
\input postpictex

\def\Z{{\rm Z\kern-.32em Z}}
\def\N{{\rm I\kern-.20em N}}
\def\R{{\rm I\kern-.17em R}}
\def\I{{\rm 1\kern-.30em I}}

\newtheorem{proposition}{Proposition}

\begin{document}
\centerline{\large \bf The Taylor property in non-negative bilinear
models}

\

\centerline{E. Gon\c{c}alves$^1$, C.M. Martins$^2$, N. Mendes-Lopes$^1$}

\

{\small  $^1${\it CMUC, Dep. of Mathematics, University of Coimbra,
Portugal}

$^2${\it Dep. of Mathematics, University of Coimbra, Portugal}}

\

\

\noindent{\bf Abstract.} The aim of this paper is to discuss the
presence of the Taylor property in the class of non-negative simple
bilinear models. Considering strictly and weakly stationary models,
we deduce autocorrelations of the process and of the square process
and analyze the presence of  the Taylor property considering several
error process distributions. The relationship between  the Taylor
property and leptokurtosis of the corresponding bilinear process is
discussed.

With the goal of extending this research to real valued bilinear
models, a simulation study is developed in a class of such models
with symmetrical innovations.

\

\noindent{\bf  Keywords.} Bilinear models; nonlinear time series;
stationarity; Taylor property.

\

\noindent{\it AMS Classification}: 62M10

\section{Introduction}
\label{Sec1}

The search for non-trivial empirical regularities in time series,
usually called stylized facts, has been the subject of several
studies in order to identify classes of time series models that
conveniently capture such empirical properties. A stylized fact
detected by  Taylor (\cite{Taylor}) when he analyzed 40 returns
series is known as  the Taylor effect. He observed that, for most of
the returns series, denoted by $X_t$ for instant $t$, the sample
autocorrelations of the absolute returns, ${\hat
\rho}_{|X|}(n)={\widehat{corr}}(|X_t|,|X_{t-n}|)$, were larger than
those of the squared returns, ${\hat
\rho}_{X^2}(n)={\widehat{corr}}(X_t^2,X_{t-n}^2)$, for $n\in
\{1,\ldots,30 \}$.

We point out that there is still little research on the theoretical
counterpart on this empirical property  due to the difficulty of
handling the true autocorrelations of time series models. For
example, this theoretical counterpart was studied by  He and
Ter\"{a}svirta (\cite{HeTer}) on conditionally Gaussian absolute value
generalized ARCH (AVGARCH) models, assuring its presence for some of
these models. More precisely, they called the theoretical relation
${\rho}_{|X|}(n)>{\rho}_{X^2}(n)$, $n\geq 1$,  the Taylor property
and concentrated their study on the autocorrelation of lag 1. More
recently, Gon\c{c}alves, Leite and Mendes-Lopes (\cite{GLML}) studied
the presence of  the Taylor property in TARCH models, concluding
that this property is satisfied when $n=1$, for some first-order
models. Generalizing these papers, Haas (\cite{Haas}) proposed a
methodology for identifying  the Taylor property in AVGARCH$(1,1)$
models at all lags.

Bilinear processes have also been proven to be suitable in financial
and physical time series modeling, namely those presenting  the
Taylor effect. Therefore, it is obviously advisable to analyze the
presence of  the Taylor property in these processes. In this paper
we consider the simple bilinear diagonal model
\begin{equation}
\label{eq1}
X_t=\beta
X_{t-k}\varepsilon_{t-k}+\varepsilon_{t},\;\;\;k>0,
\end{equation}
where $\beta$ is a real parameter and $(\varepsilon_t, t\in \Z)$ an
error process. We state sufficient conditions for the strict and
weak stationarity of the processes $X=(X_t,t\in\Z)$ and
$X^2=(X_t^2,t\in\Z)$, and we derive expressions for the moments of
$X$ up to the 4th order.

When dealing with bilinear models it is common to assume that
$\varepsilon_t$, $t\in \Z$, are normally distributed. However, there
has been considerable interest in non-negative time series models.
For instance, Pereira and Scotto (\cite{PS}) studied some properties
of the simple first-order bilinear diagonal
 model ($k=1$) driven by exponentially distributed innovations.

In this paper, we analyze the presence of  the Taylor property  when
$n=1$ in the non-negative first-order bilinear time series model
considering several distributions for the error process, which are
chosen according to the kurtosis  value as we have observed that the
Taylor property is related with the value of this parameter.

Based on a simulation study, we also analyze the presence of the
Taylor property in the class of real valued first-order bilinear
diagonal models with symmetrical innovations.

\section{Stationarity of $X$ and $X^2$}
\label{Sec2}

In this section we consider the simple bilinear model defined by
(\ref{eq1}) where $(\varepsilon_{t},t\in \Z)$ is a sequence of
i.i.d. random variables.  Let $\mu_i=E(\varepsilon_t^i)$, $i\in \N$.

\begin{proposition}
\label{Propos1} Suppose that $\mu_4$ and $E(\ln |\varepsilon_t|)$
exist. If $\beta^2\, \mu_2<1$ then the process $X$ is strictly and
weakly stationary.
\end{proposition}

\

\noindent {\bf Proof.} To prove the strict stationarity of process
$X$, we start by proving that $X_t=Y_t$, a.s., with
$$Y_t=\varepsilon_t+\sum_{n=1}^{+\infty}T_n,$$ where, for each $n\in \N$, $T_n=T_n(t)$ is given
by
$$T_n=\beta^n
\varepsilon_{t-nk}\prod_{j=1}^{n}\,\varepsilon_{t-jk}.$$

Let us begin by verifying that the series $\sum_{n=1}^{\infty}T_n$
is a.s. convergent. Using the ergodic theorem, we can assure that
the limit $\,{\displaystyle \lim_{n\rightarrow
+\infty}}\frac{1}{n}\ln \left|\beta^n
\prod_{j=1}^{n}\varepsilon_{t-jk}  \right|$ exists and that
${\displaystyle \lim_{n\rightarrow +\infty}}\frac{1}{n}\ln
\left|\beta^n \prod_{j=1}^{n}\varepsilon_{t-jk}  \right|=\ln |\beta|
+ E(\ln |\varepsilon_t |)$.

We can observe that $\frac{1}{n}\ln |T_n|= \frac{1}{n}\ln
\left|\beta^n \prod_{j=1}^{n}\varepsilon_{t-jk}
\right|+\frac{1}{n}\ln |\varepsilon_{t-nk}|$.  Since \linebreak
${\displaystyle \lim_{n\rightarrow +\infty}(a.s.)\frac{1}{n}\ln
|\varepsilon_{t-nk}|=0}$, we have
$${\displaystyle \lim_{n\rightarrow +\infty}(a.s.)\frac{1}{n}\ln
|T_n|=\ln |\beta| + E(\ln |\varepsilon_t |)}.$$

On the other hand, the condition $\beta^2\mu_2<1$ implies $2\ln
|\beta |<-\ln E(\varepsilon_t^2)$. Applying Jensen's inequality to
the random variable $\varepsilon_t^2$ and taking into account that
$E(|\ln |\varepsilon_t||)<+\infty$, we obtain
$\gamma=\ln|\beta|+E(\ln |\varepsilon_t|)<0$.

Consequently
$$\lim_{n\rightarrow \infty}(a.s.){|T_n(t)|}^{1/n}=
\exp\gamma<1,$$ which implies that the series
$\sum_{n=1}^{\infty}T_n$ is a.s. convergent; so $(Y_t,t\in\Z)$ is a
strictly stationary process, as it is a measurable function of the
independent random variables $\varepsilon_s$, $s\leq t$. Moreover,
it is easy to verify that the process $(Y_t, t\in\Z)$ satisfies
Equation (\ref{eq1}).

This solution is the unique strictly stationary solution of
(\ref{eq1}). In fact, using (\ref{eq1}) recursively, we obtain $$
X_t= \varepsilon_t+\sum_{i=1}^{n}T_i+\beta^{n+1}X_{t-(n+1)k}
\prod_{j=0}^{n}\varepsilon_{t-(j+1)k}, \,\; n=0,1,\ldots
$$
\noindent with ${\sum_{n=1}^{0}T_n=0}$, for each $t\in\Z$, and
taking limits, any strictly stationary solution of (\ref{eq1})
satisfies
$$X_t=Y_t+\lim_{n\rightarrow
+\infty}(a.s.)  \beta^{n+1} X_{t-(n+1)k}
\prod_{j=0}^{n}\varepsilon_{t-(j+1)k}.$$

Let ${ Z_n(t)=\beta^n
X_{t-nk}\prod_{j=0}^{n-1}\varepsilon_{t-(j+1)k}}$. It is easy to
verify that $${\displaystyle  \lim_{n\rightarrow+\infty}}
(a.s.)\frac{1}{n}\ln |Z_n(t)|=\gamma<0.$$ Then
$$ \lim_{n\rightarrow+\infty}(a.s.)
|Z_n(t)|=\lim_{n\rightarrow+\infty}(a.s.)\exp
\left[n\left(\frac{1}{n}\ln |Z_n(t)| \right)\right]=0,$$ \noindent
which implies ${\displaystyle \lim_{n\rightarrow+\infty}(a.s.)
Z_n(t)=0}$. So,   $(X_t, t\in\Z)$ is  strictly stationary, as
$X_t=Y_t$, a.s..

To prove the weak stationarity, we now verify that
$\,E(Y_t^2)<+\infty$. We have
\begin{eqnarray}
E\left(Y_t^2\right) & = & E\left[\left( \varepsilon_t+\sum_{i=1}^{+\infty}T_i \right)^2\right] \nonumber \\
& \leq &
E\left(\varepsilon_t^2\right)+2\sum_{i=1}^{\infty}E\left(|\varepsilon_t|
|T_i|\right)+\sum_{i=1}^{\infty}
\sum_{j=1}^{\infty}E\left(|T_iT_j|\right). \label{eq2}
\end{eqnarray}
Under the given conditions, each series in (\ref{eq2}) is
convergent. In fact, let us consider, for example, the series
$\,{\displaystyle
\sum_{i=1}^{\infty}\sum_{j=1}^{\infty}E\left(|T_iT_j|\right)}$.

For each  $i,j\in \N$, we have
\begin{eqnarray*} E\left(|T_iT_j|\right)  & \leq &
|\beta|^{i+j}\left[E\left(\varepsilon_{t-ik}^4\varepsilon_{t-k}^2
\varepsilon_{t-2k}^2\ldots
\varepsilon_{t-(i-1)k}^2\right)\right]^{1/2} \\ \hspace*{1.0cm} & &
\left[E\left(\varepsilon_{t-jk}^4\varepsilon_{t-k}^2
\varepsilon_{t-2k}^2\ldots
\varepsilon_{t-(j-1)k}^2\right)\right]^{1/2} \\ & = &
\mu_4\mu_2^{-1}\left[\left( \beta^2\mu_2\right)^{1/2}\right]^{i+j},
\end{eqnarray*}
using Schwarz's inequality and the independence of the r.v.'s
$\varepsilon_t$, $t\in \Z$. As $\left( \beta^2\mu_2\right)^{1/2}<1$,
the series is convergent.

Taking into account the equality $X_t=Y_t$, a.s., and the strict
stationarity of the process $X$, we conclude that $E(X_t^2)$ exists
and that $X$ is weakly stationary.

\

\begin{proposition}
\label{Propos2}Suppose that $E(\ln |\varepsilon_t|)$ and  $\mu_8$
exist. If  $\beta^4\, \mu_4<1$ then the process $X^2$ is strictly
and weakly stationary.
\end{proposition}

\

\noindent {\bf Proof.} The condition $\beta^4\mu_4<1$ implies
$\beta^2\mu_2<1$, using Schwarz's inequality, which implies the
strict stationarity of $X$ and, consequently, of $X^2$. The proof of
the  weak stationarity of $X^2$ is analogous to the previous one. We
have
\begin{eqnarray*} & & E\left( Y_t^4 \right) \leq
E\left(\varepsilon_t^4\right)
+\sum_{i=1}^{\infty}\sum_{j=1}^{\infty}\sum_{p=1}^{\infty}
\sum_{q=1}^{\infty}E\left(\left|T_iT_jT_pT_q\right|\right) +4
\sum_{i=1}^{\infty}E\left(|\varepsilon_t^3|\left|T_i\right|\right)
\\ & & \hspace*{1.5cm}
+4\sum_{i=1}^{\infty}\sum_{j=1}^{\infty}\sum_{p=1}^{\infty}
E\left(|\varepsilon_t|\left|T_iT_jT_p\right|\right)+6
\sum_{i=1}^{\infty}\sum_{j=1}^{\infty}E\left( \varepsilon_t^2
\left|T_iT_j\right|\right).
\end{eqnarray*}

Let us consider, for example, the series ${\displaystyle
\,\sum_{i=1}^{\infty}\sum_{j=1}^{\infty}\sum_{p=1}^{\infty}
\sum_{q=1}^{\infty}E\left(\left|T_iT_jT_pT_q\right|\right)} $, which
is a sum of series of the types
\begin{enumerate}
\item[({\it i})] ${\displaystyle \sum_{i=1}^{\infty}
\sum^{\infty}_{j=i+1} \sum_{p=1}^{\infty}
\sum^{\infty}_{q=p+1}E\left(\left|T_iT_jT_pT_q\right|\right)}$
\item[({\it ii})] ${\displaystyle \sum_{i=1}^{\infty}
\sum_{p=1}^{\infty}E\left( T_i^2T_p^2 \right)}$
\item[({\it iii})] ${\displaystyle \sum_{i=1}^{\infty}
\sum^{\infty}_{j=i+1} \sum_{p=1}^{\infty}
E\left(\left|T_iT_j\right|T_p^2\right)}$.
\end{enumerate}

Concerning ({\it i}), as  $j>i$ and $q>p$, we have
\begin{eqnarray*}
& &  E\left(\left|T_iT_jT_pT_q\right|\right) =
 E\left[\left(|T_iT_j|\right)\left(|T_pT_q\right)\right]\\
 & & \hspace*{0.5cm}\leq
\left[E\left(\varepsilon_{t-ik}^4\varepsilon_{t-l}^4
\varepsilon_{t-l-k}^4\ldots
\varepsilon_{t-l-(i-1)k}^4\varepsilon_{t-jk}^2
\varepsilon_{t-l-ik}^2\ldots\varepsilon_{t-l-(j-1)k}^2
\right)\right]^{1/2}
\\ & & \hspace*{1.25cm}\left[
E\left(\varepsilon_{t-pk}^4\varepsilon_{t-l}^4
\varepsilon_{t-l-k}^4\ldots
\varepsilon_{t-l-(p-1)k}^4\varepsilon_{t-qk}^2
\varepsilon_{t-l-pk}^2\ldots\varepsilon_{t-l-(q-1)k}^2
\right)\right]^{1/2},
\end{eqnarray*}
\noindent using Schwarz's inequality.

Taking into account the independence of the random variables
$\varepsilon_t$, we have, for $i,j\in\N$, $ \,j>i$,
$$E\left(\varepsilon_{t-ik}^4\,\varepsilon_{t-l}^4\,
\varepsilon_{t-l-k}^4\ldots
\varepsilon_{t-l-(i-1)k}^4\,\varepsilon_{t-jk}^2\,
\varepsilon_{t-l-ik}^2\ldots\varepsilon_{t-l-(j-1)k}^2
\right)=\mu_4^{i+1}\mu_2^{j-i+1}.$$

\noindent Then \begin{eqnarray*} & & \sum_{i=1}^{\infty}
\sum^{\infty}_{j=i+1} \sum_{p=1}^{\infty} \sum^{\infty}_{
q=p+1}E\left(\left|T_iT_jT_pT_q\right|\right)
\\ & & \hspace*{0.8cm} \leq \sum_{i=1}^{\infty}
\sum^{\infty}_{j=i+1} \sum_{p=1}^{\infty} \sum^{\infty}_{
q=p+1}{|\beta|}^{i+j+p+q}\left(\mu_4^{i+p+2}\mu_2^{j-i+q-p+2}
\right)^{1/2}  \\ & & \hspace*{1.4cm} =\sum_{i=1}^{\infty}
\sum^{\infty}_{j=i+1} \sum_{p=1}^{\infty}
\sum^{\infty}_{q=p+1}\mu_2\mu_4\left[\left(\beta^4\mu_4\right)^{1/2}
\right]^{i+p} \left[\left(\beta^2\mu_2\right)^{1/2}
\right]^{[(j+q)-(i+p)]}.
\end{eqnarray*}
\noindent As $\,(\beta^4\mu_4)^{1/2}<1$ and
$\,(\beta^2\mu_2)^{1/2}<1$, the series in ({\it i}) is convergent.
The convergence of the series ({\it ii}) and ({\it iii}) is proved
in a similar way. Then we conclude that $\,E(X_t^4)<+\infty$,
$\,t\in\Z$. As the process $X^2$ is strictly stationary and
$E(X_t^4)$ exists, then it is weakly stationary.

\

\section{Moments up to the 4th order}
\label{Sec3}

Under the same conditions of Section \ref{Sec2}, we now evaluate the
moments up to the 4th order of the process $X$ given by (\ref{eq1})
where $(\varepsilon_{t},t\in \Z)$ is a sequence of i.i.d. random
variables, and $\mu_i=E\{\varepsilon_t^i\}$, $i\in \N$.

\begin{proposition}
If $\beta^4 \mu_4<1$ and $\mu_8$ exists then the $n$th moment of
$X_t$, $n\leq 4$, can be expressed as
$$E(X_t^n)=\sum_{i=0}^{n}{n\choose
i}\beta^{n-i}\,\mu_i\,E(X_t^{n-i}\varepsilon_t^{n-i}),$$ \noindent
where $$\displaystyle E(X_t^n\varepsilon_t^n) = \frac{1}{1-\beta^n
\mu_n}\; \sum_{i=1}^{n}{n\choose
i}\beta^{n-i}\,\mu_{n+i}\,E(X_t^{n-i}\varepsilon_t^{n-i}),\,\;\;n\leq
4.$$
\end{proposition}

\

\noindent {\bf Proof.} For $n\leq 4$, we have
\begin{eqnarray*}E(X_t^n)& = &
\sum_{i=0}^{n}{n\choose
i}\beta^{n-i}\,E\left[\varepsilon_t^{i}\left(X_{t-k}\varepsilon_{t-k}\right)^{n-i}\right] \\
& = & \sum_{i=0}^{n}{n\choose
i}\beta^{n-i}\,\mu_i\,E(X_t^{n-i}\varepsilon_t^{n-i}),
\end{eqnarray*}
\noindent since the process  $(X_t\varepsilon_t, t\in\Z)$ is
strictly stationary due to the fact that $X_t\varepsilon_t$ is a
measurable function of $\varepsilon_{t}, \varepsilon_{t-1}, \ldots$.
Now we need to evaluate $E(X_t^n\varepsilon_t^n)$, $n\leq 4$.
\begin{eqnarray*}E(X_t^n\varepsilon_t^n) & = & \sum_{i=0}^{n}{n\choose
i}\beta^{n-i}\,E\left[\varepsilon_t^{i}\left(X_{t-k}\varepsilon_{t-k}\right)^{n-i}\varepsilon_t^{n}\right] \\
& = & \sum_{i=0}^{n}{n\choose
i}\beta^{n-i}\,E\left(\varepsilon_t^{n+i}\right)\,E
\left(X_{t}^{n-i}\varepsilon_{t}^{n-i}\right) \\
& = & \beta^{n}\,\mu_n\,E
\left(X_{t}^{n}\varepsilon_{t}^{n}\right) +\sum_{i=1}^{n}{n\choose
i}\beta^{n-i}\,\mu_{n+i}\,E
\left(X_{t}^{n-i}\varepsilon_{t}^{n-i}\right).
\end{eqnarray*}

Then $$\displaystyle E(X_t^n\varepsilon_t^n) = \frac{1}{1-\beta^n
\mu_n}\; \sum_{i=1}^{n}{n\choose
i}\beta^{n-i}\,\mu_{n+i}\,E(X_t^{n-i}\varepsilon_t^{n-i}).$$

\

It is easy to verify that $\,E(X_t\varepsilon_t)=\mu_2/(1-\beta
\mu_1)$. Recursively, we obtain $E(X_t^n\varepsilon_t^n)$,
$n=1,2,3$, and, finally, we achieve $E(X_t^n)$,  $n\leq 4$.

We note that $\beta^4\mu_4<1$ implies $|\beta^n\mu_n|<1$, $n=1,2,3$,
using Schwarz's inequality.

\section{The Taylor property in first-order non-negative bilinear
models} \label{Sec4}

In this section we consider the first-order non-negative bilinear
model
\begin{equation}
X_t=\beta X_{t-1}\varepsilon_{t-1}+\varepsilon_{t}, \;\;\;t\in
\Z,
\label{eq3}
\end{equation}
\noindent where  $\beta>0$ and $(\varepsilon_{t},t\in \Z)$ is a
sequence of non-negative i.i.d. random variables.

We assume that $E(\ln \varepsilon_t)$ and $\mu_8$ exist and that
$\beta^4\, \mu_4<1$ in order to guarantee that both processes, $X$
and $X^2$, are strictly and weakly stationary.

In this context,  the Taylor property for $n=1$ establishes that
$\rho_{X}(1)>\rho_{X^2}(1)$, where $\rho_{X}(1)$ and $\rho_{X^2}(1)$
denote, respectively, the autocorrelations of lag 1 of the processes
$X$ and $X^2$. In order to obtain these autocorrelations, it is
enough to evaluate $E(X_tX_{t-1})$ and $E(X_t^2X_{t-1}^2)$ since we
derived $E(X_t^i)$, $i=1,2,3,4$, in the previous section.
 Using (\ref{eq3}) and the stationarity of the
involved processes, we have
\begin{eqnarray*}
E(X_tX_{t-1}) & = & \beta
E(X_{t}^2\varepsilon_{t})+E(X_{t-1}\varepsilon_{t}) \\
& = &  \beta
E(\beta^2X_{t-1}^2\varepsilon_{t-1}^2\varepsilon_{t}+2\beta
X_{t-1}\varepsilon_{t-1}\varepsilon_{t}^2+\varepsilon_{t}^3)+E(X_{t-1}\varepsilon_{t}).
\end{eqnarray*}
Taking into account the independence of the random variables
$\varepsilon_{t}$, $t\in\Z$, and the strict stationarity of the
related processes, we have
$E(X_{t-1}^2\varepsilon_{t-1}^2\varepsilon_{t})=\mu_1
E(X_{t}^2\varepsilon_{t}^2)$ and
$E(X_{t-1}\varepsilon_{t-1}\varepsilon_{t}^2)=\mu_2
E(X_{t}\varepsilon_{t})$. Then
$$E(X_tX_{t-1})=\beta^3 \mu_1  E(X_{t}^2\varepsilon_{t}^2)+2\beta^2\mu_2E(X_{t}\varepsilon_{t})
+\mu_1 E(X_t)+\beta\mu_3.$$ \indent Using an analogous procedure, we obtain
\begin{eqnarray*}
E(X_t^2X_{t-1}^2)
& = &  \beta^4 E_1+ 2\beta^3 E_2 + 2\beta^3\mu_1 E_3+4\beta^2 \mu_1 E_4 +\beta^2 E_5+2\beta \mu_1 E_6  \\
& & \hspace*{1.6cm}+ \beta^2 \mu_2
E(X_t^2\varepsilon_t^2)+2\beta\mu_1 \mu_2
E(X_t\varepsilon_t)+\mu_2^2,
\end{eqnarray*}
\noindent where
\begin{eqnarray*}
E_1 & = & E(X_t^2X_{t-1}^2\varepsilon_t^2\varepsilon_{t-1}^2) \;=\; \beta^2 \mu_2 E(X_t^4\varepsilon_t^4)+2\beta \mu_3 E(X_t^3\varepsilon_t^3)
+\mu_4E(X_t^2\varepsilon_t^2) \\
E_2 & = & E(X_t^2X_{t-1}\varepsilon_t^3\varepsilon_{t-1})  \;=\;  \beta^2 \mu_3 E(X_t^3\varepsilon_t^3)
+2\beta \mu_4 E(X_t^2\varepsilon_t^2)+\mu_5 E(X_t\varepsilon_t) \\
E_3 & = & E(X_tX_{t-1}^2\varepsilon_t\varepsilon_{t-1}^2) \;=\;  \beta \mu_1 E(X_t^3\varepsilon_t^3)
+\mu_2E(X_t^2\varepsilon_t^2) \\
E_4 & = & E(X_tX_{t-1}\varepsilon_t^2\varepsilon_{t-1})  \;=\; \beta
\mu_2 E(X_t^2\varepsilon_t^2) +\mu_3 E(X_t\varepsilon_t) \\
E_5 & = &
E(X_t^2\varepsilon_t^4)  \;=\;  \beta^2 \mu_4
E(X_t^2\varepsilon_t^2)
+2\beta\mu_5 E(X_t\varepsilon_t)+\mu_6 \\
E_6 & = & E(X_t\varepsilon_t^3)  \;=\;  \beta\mu_3 E(X_t\varepsilon_t)+\mu_4.
\end{eqnarray*}
Finally, the results of the previous section allow us to obtain the
values of $E(X_tX_{t-1})$ and $E(X_t^2X_{t-1}^2)$ in terms of the
moments of $\varepsilon_t$.

In the following, we investigate the presence of  the Taylor
property in Model (\ref{eq3}), considering some non-negative
distributions for the error process, namely, the uniform
distribution in $]0,\alpha[$, the exponential distribution in
$]0,+\infty[$ with mean $\alpha$, and the Pareto distribution with
density $\displaystyle f(x)=\frac{\nu \alpha^{\nu
}}{x^{\nu+1}}\I_{]\alpha,+\infty[}(x)$, for $\nu=12$ and $\nu=9$. In
all cases, $\alpha$ is a non-negative parameter and the condition
$E(|\ln \varepsilon_t|)<+\infty$ is satisfied.

The choice of these distributions takes into account the fact that
 the Taylor property seems to be related with the kurtosis value of
the process. In this sense, we choose four distributions with
significantly different behavior as regards their tails. We point
out that the uniform and exponential distributions have constant
kurtosis values, while the kurtosis of the Pareto distribution
depends on the value of the parameter $\nu$. Consequently, valid
comparisons may be made separately between the first two
distributions, uniform and exponential, and then between the two
referred Pareto distributions.

We also point out that, in all cases, the condition $\beta^4 \,
\mu_4<1 $ and the values of $\rho_{X}(1)$ and $\rho_{X^2}(1)$ can be
written in terms of $r=\alpha \beta$.

In each case, we also present the value of the kurtosis of the
process $X$ given by (4.1), which also depends on $r=\alpha \beta$,
as well as the corresponding graphic representation as a function of
$r$.

\vspace{0.2cm} \noindent{\bf Error process with uniform distribution
in $]0,\alpha[$}

\noindent In this case, the condition $\beta^4 \, \mu_4<1$ is
equivalent to $0<r<\sqrt[4]{5}\simeq  1.495$ and we obtain
\begin{eqnarray*}
\rho_X(1) & = & \frac{r(-180 +
120 r - 51 r^2 - 4 r^3 + r^4)}{-180 + 180 r - 177 r^2 +
 12 r^3 + 7 r^4}\\
\rho_{X^2}(1)  & = & -\frac{r}{12}\,\frac{N_u(r)}{D_u(r)},
\end{eqnarray*} with \begin{eqnarray*} N_u(r) & = & -604800 - 480600
r - 155700 r^2 - 257400 r^3 - 2490 r^4 +
     48525 r^5 \\  & \; & \hspace*{0.5cm} - 6270 r^6
    + 6810 r^7 + 10620 r^8 + 11384 r^9 +
     4012 r^{10} - 586 r^{11} \\  & \; & \hspace*{1.0cm}+ 94 r^{12} - 53 r^{13} + 6 r^{14}\\
D_u(r)& = & 50400 + 12600 r + 35700 r^2 + 40200 r^3 + 13490 r^4 +
14015 r^5 +
   8360 r^6\\
     & \; & \hspace*{0.25cm}- 5210 r^7 - 5999 r^8 - 2407 r^9 - 720 r^{10} + 114 r^{11} +
   177 r^{12} - 8 r^{13}. \end{eqnarray*}

\vspace*{5mm}
\begin{figure}[h]
\hspace{0.8cm}
\epsfxsize=2.2in \epsffile{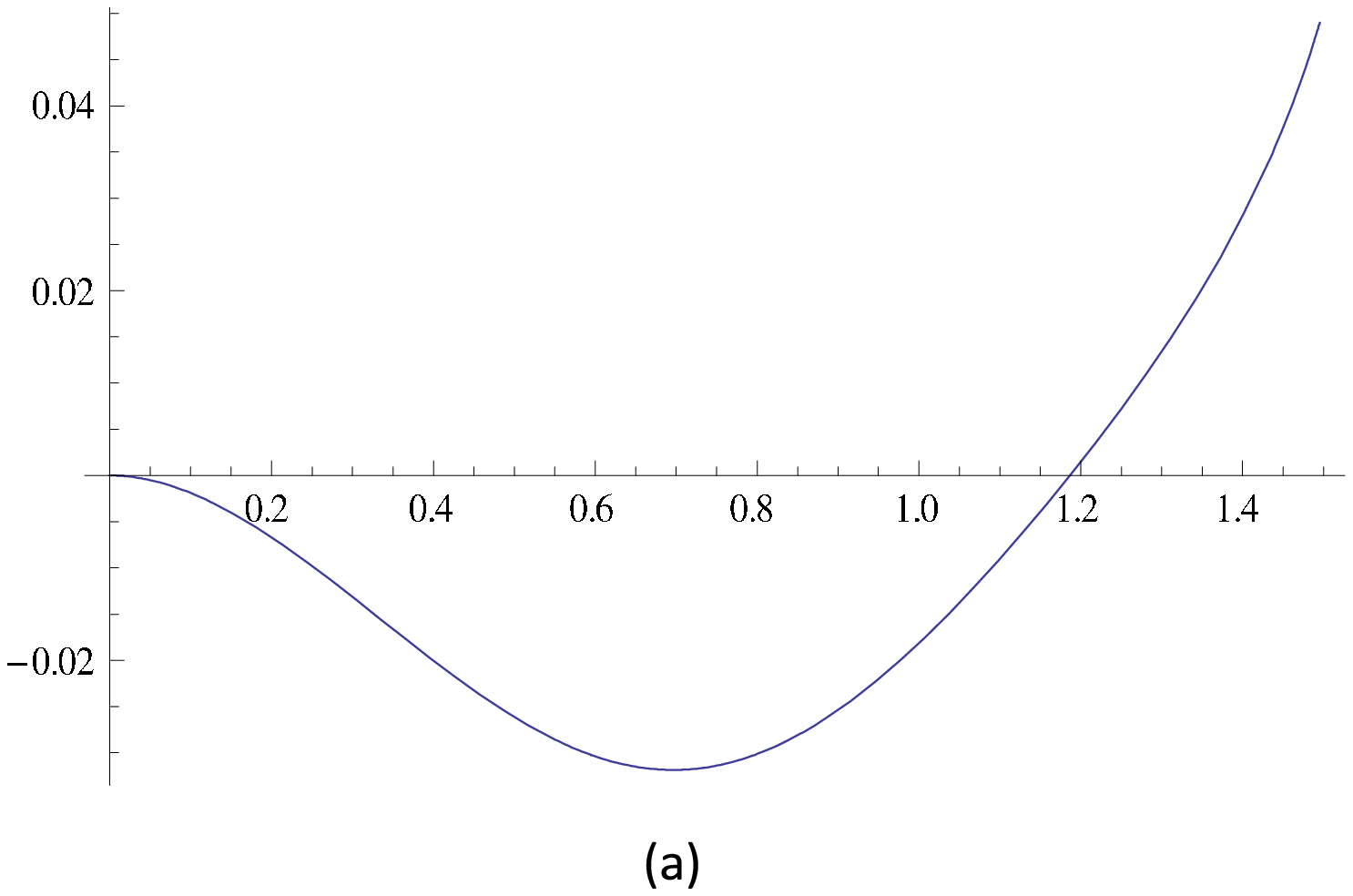}
\hspace{0.5cm}
\epsfxsize=2.2in \epsffile{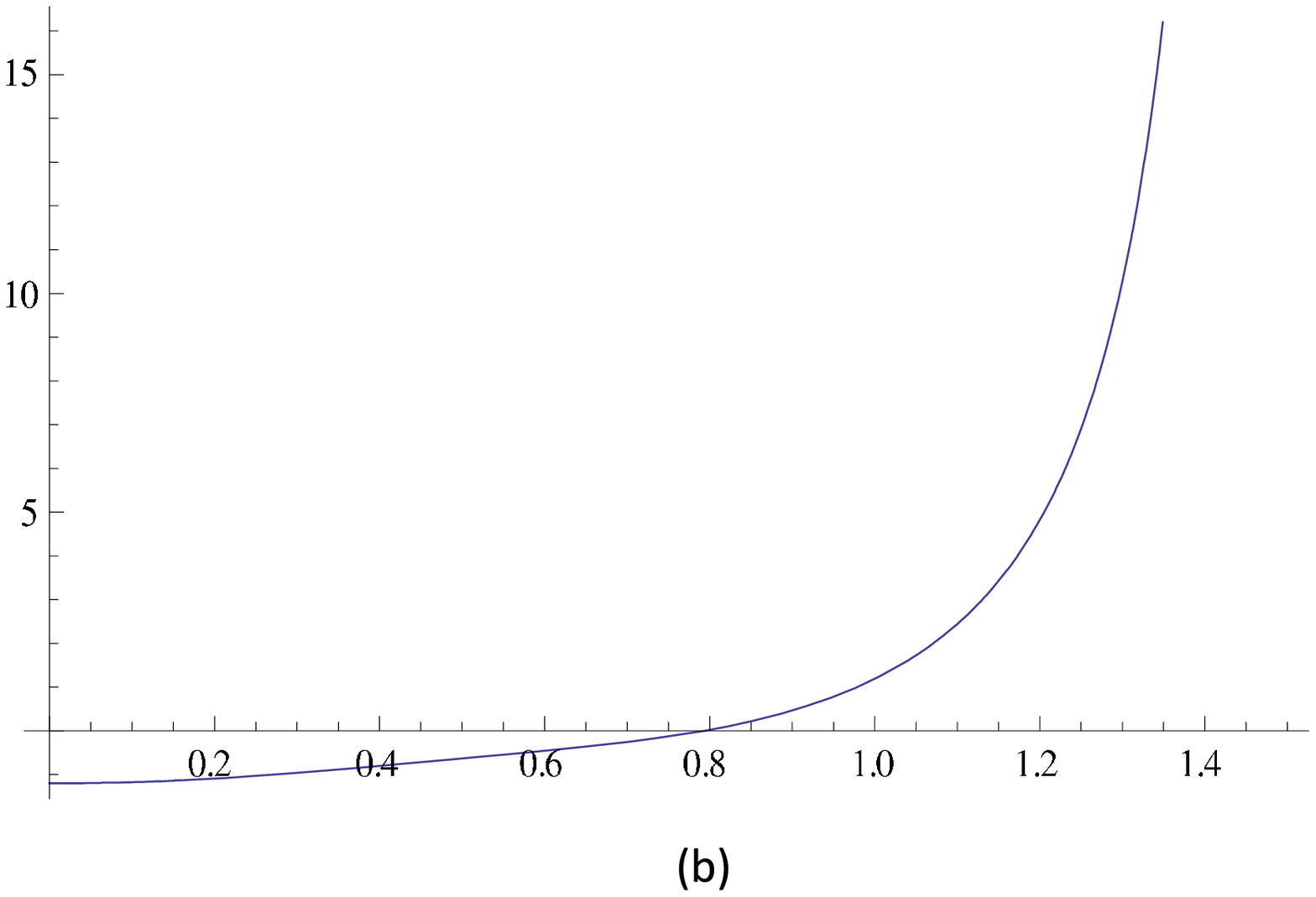}
\caption{Graphs from $\rho_X(1)-\rho_{X^2}(1)$
{\sf (a)} and $K_u(r)$ {\sf (b)},  with $0<r<\sqrt[4]{5}$
 \label{fig1}}
\end{figure}
\vspace*{5mm}

From Figure~\ref{fig1}(a), we can see that  the Taylor property is
present for values of $r$ in the interval $]1.1868987,
\sqrt[4]{5}[$. So, for a fixed $\alpha$,  the Taylor property is
achieved for parameterizations of Model (\ref{eq3}) such that
$$\beta\in
\left]\frac{1.1868987}{\alpha},\frac{\sqrt[4]{5}}{\alpha}\right[,$$
\noindent where the value $1.1868987$ was obtained with an
approximation error inferior to $5\times 10^{-9}$.

For Model (\ref{eq3}) with such an error process, the kurtosis is
given by
$$K_u(r)=\frac{-3 (-3 + r^2)}{7 (-4 + r^3) (-5 + r^4)}\frac{N_u^*(r)}{D_u^*(r)}-3,$$ where
\begin{eqnarray*}
N_u^*(r) & = & 907200 - 1814400 r + 4284000 r^2 - 4510800 r^3 +
      3254460 r^4 \\ & \; & \hspace*{0.5cm} - 2030520 r^5
+ 1973540 r^6 - 617175 r^7 -
      185700 r^8 + 371005 r^9 \\  & \; & \hspace*{1.5cm} - 236308 r^{10}
 + 78747 r^{11} -
      11496 r^{12} + 511 r^{13} \\
D_u^*(r) & = &  (-180 +
     180 r - 177 r^2 + 12 r^3 + 7 r^4)^2.
\end{eqnarray*}
From Figure~\ref{fig1}(b), we observe that the kurtosis of this
model is an increasing function of $r$ and, for large values of the
kurtosis,  the Taylor property occurs.

\vspace*{0.2cm} \noindent{\bf Error process with  exponential
distribution with mean $\alpha$ (in $]0,+\infty[$)}

\begin{figure}[h]
\hspace{0.8cm}
\epsfxsize=2.2in \epsffile{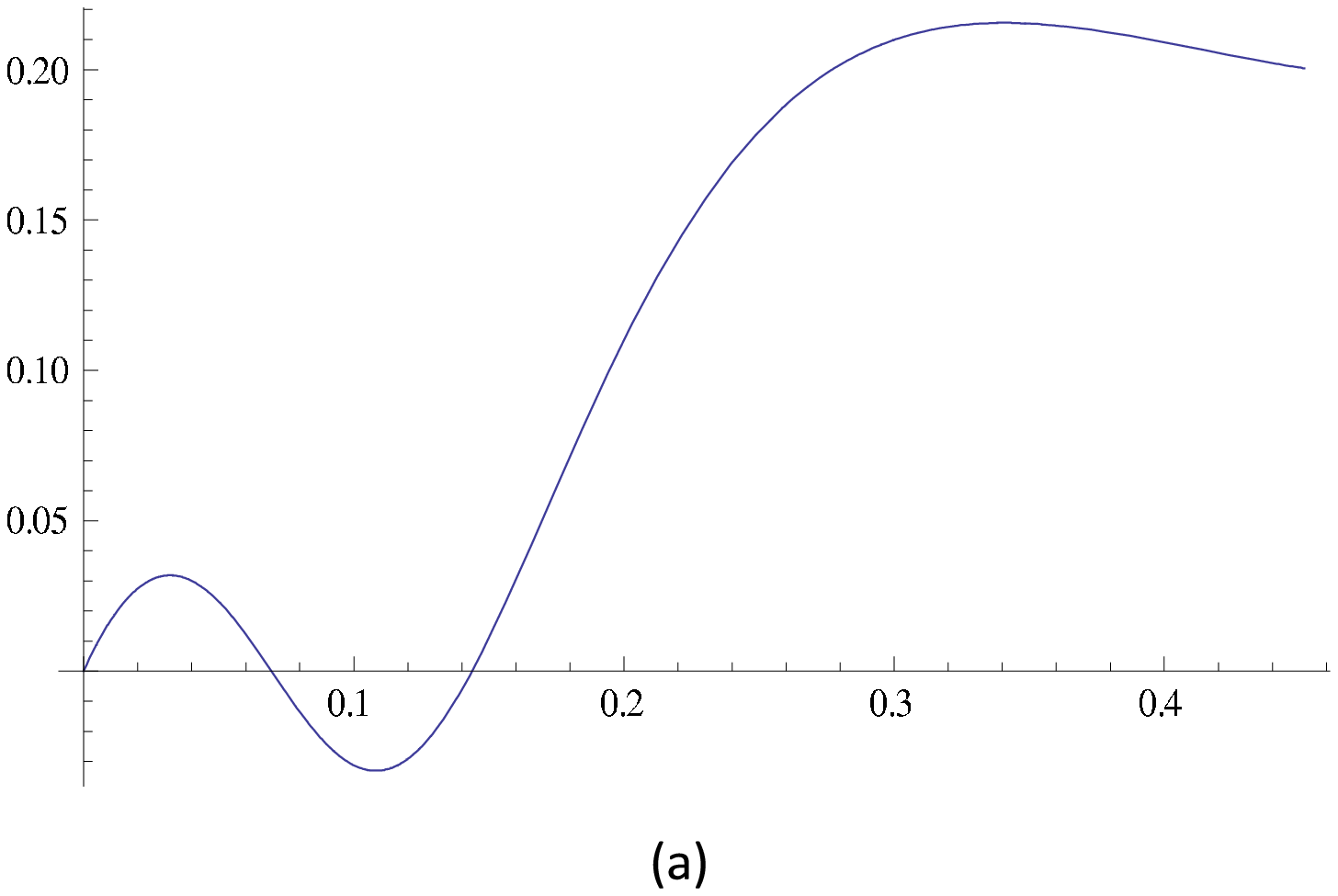}
\hspace{0.2cm}
\epsfxsize=2.2in \epsffile{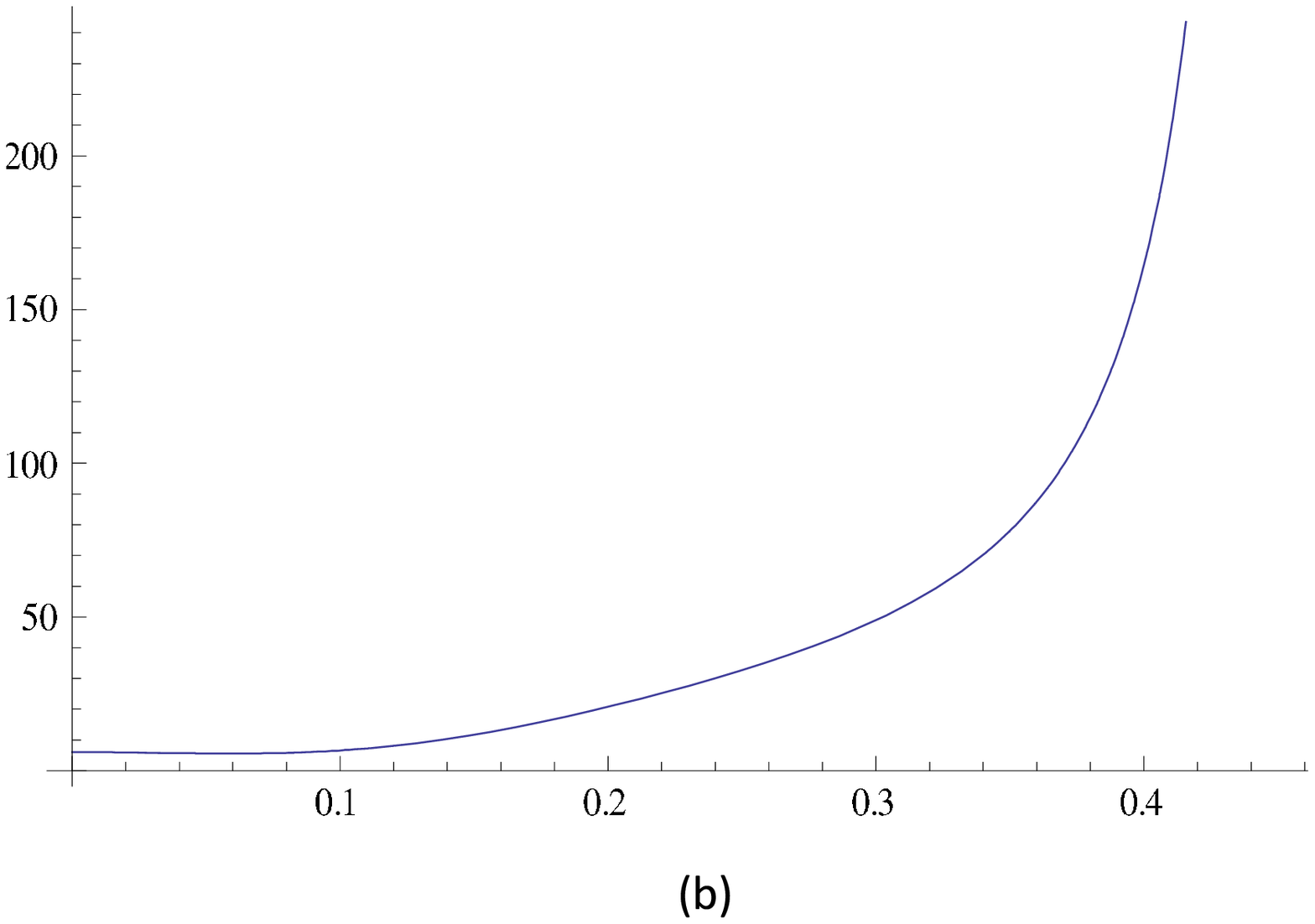}
\caption{Graphs from $\rho_X(1)-\rho_{X^2}(1)$
{\sf (a)} and $K_e(r)$  {\sf (b)}, with
$0<r<\frac{1}{\sqrt[4]{24}}$\label{fig2}}
\end{figure}

\noindent The condition $\beta^4 \, \mu_4<1$ is now equivalent to
$0<r<\frac{1}{\sqrt[4]{24}}\simeq 0.4518$. In this case,
\begin{eqnarray*}
\rho_X(1) & = & \frac{2r(
        2-3r+7r^2-6r^3+2r^4)}{1-2r+19r^2-20r^3+6r^4}\\
\rho_{X^2}(1)  & = & 2r\,\frac{N_e(r)}{D_e(r)}.
\end{eqnarray*}
\noindent with
\begin{eqnarray*} N_e(r) & = &  -5 - 80 r + 65 r^2 - 112 r^3 - 1184 r^4 - 5774 r^5 +
10848 r^6 +12720 r^7 \\ & \; & \hspace*{0.5cm} - 9408 r^8
 - 17880 r^9 - 16272 r^{10} + 52992 r^{11} +
 9216 r^{12} \\  & \; & \hspace*{1.0cm} - 46656 r^{13} + 17280 r^{14} \\
D_e(r)& = &
-5 + 2 r - 21 r^2 - 602 r^3 - 9060 r^4 + 11126 r^5 + 13252 r^6-
 26448 r^7\\ & \; & \hspace*{0.3cm}+ 16368 r^8
 + 13896 r^9 - 12192 r^{10} + 13824 r^{11} -
 12672 r^{12} + 4032 r^{13}.
\end{eqnarray*}
So, when the errors are exponentially distributed with mean
$\alpha$, Model (\ref{eq3}) presents  the Taylor property for
parameterizations such that
$$\beta \in \left]0,\frac{0.0695566}{\alpha}\right[\; \cup\;
\left]\frac{0.1437879}{\alpha},
\frac{1}{\sqrt[4]{24}\,\alpha}\right[,$$ where the values $0.0695566$ and
 $0.1437879$ were obtained with an
approximation error inferior to $5\times 10^{-8}$. This conclusion
is illustrated in Figure~\ref{fig2}(a). In Figure~\ref{fig2}(b), we
have the graphic representation of the kurtosis of model (\ref{eq3})
with exponential errors, which is given by
$$K_e(r)=\frac{-3 (-1 +2 r^2)}{(-1 +
     6 r^3) (-1 + 24 r^4)}\frac{N_e^*(r)}{D_e^*(r)}-3,$$ where
\begin{eqnarray*}
N_e^*(r) & = & 3 - 12 r + 52 r^2 - 134 r^3 + 11815 r^4 -
      36752 r^5 + 44802 r^6 + 1062 r^7\\
 & \; & \hspace*{0.3cm}
     - 42648 r^8   + 17028 r^9 +
      12240 r^{10} + 5616 r^{11} - 17280 r^{12} + 6048 r^{13} \\
D_e^*(r) & = &  (1 - 2 r + 19 r^2 - 20 r^3 + 6 r^4)^2.
\end{eqnarray*}

As in the previous case, the kurtosis of Model (\ref{eq3}) is an
increasing function of $r$ and large kurtosis values  correspond to
large values of the difference $\rho_X(1)-\rho_{X^2}(1)$.

We also observe that the kurtosis of the process $X$ is larger when
the errors are exponentially distributed than when the errors are
uniformly distributed, corresponding to an analogous relation
between the kurtosis of those error processes.  The Taylor property
seems to emerge in a relatively stronger way when the kurtosis of
$X$ increases.

\vspace*{0.2cm} \noindent{\bf Error process with Pareto density
$\displaystyle
f(x)=\frac{12\alpha^{12}}{x^{13}}\I_{]\alpha,+\infty[}(x)$}

\noindent  The region of existence of the autocorrelations in terms
of $r=\alpha\beta$ is now defined by $0<r
<\sqrt[4]{\frac{2}{3}}\simeq 0.9036$. We have
\begin{eqnarray*}
\rho_X(1) & = & \frac{44 r (6050 - 10230 r + 13035 r^2 - 7524 r^3 + 1296
r^4)}{3 (36300 -
   79200 r + 219255 r^2 - 171160 r^3 + 29472 r^4)}\\
\rho_{X^2}(1)  & = & \frac{r}{55}\,\frac{N_{p12}(r)}{D_{p12}(r)},\end{eqnarray*} \noindent with
\begin{eqnarray*}N_{p12}(r) & = & -7043652000 - 5638479000 r - 1900483200 r^2 - 6228372150 r^3\\
& \; & \hspace*{0.1cm} -3064649280 r^4
 + 2622844140 r^5 + 24533447400 r^6 \\& \; & \hspace*{0.2cm} +
 19854650865 r^7 + 11360213480 r^8
- 16340416020 r^9 \\& \; & \hspace*{0.3cm} -
 30235824828 r^{10} + 23037530976 r^{11} + 7650162960 r^{12}\\& \; & \hspace*{0.4cm} -
 11215587456 r^{13} + 2802615552 r^{14}\\
D_{p12}(r) & = & -58697100 + 14229600 r - 142425360 r^2 - 468153840 r^3 \\ & \; & \hspace*{0.1cm}-
 218936564 r^4 + 536116224 r^5  + 616017864 r^6 \\& \; & \hspace*{0.2cm} + 374454192 r^7 +
 130906149 r^8 - 805701976 r^9 \\ & \; & \hspace*{0.4cm}- 15605040 r^{10}
 + 401099652 r^{11} \\ & \; & \hspace*{0.4cm} -
 245871648 r^{12} + 48736320 r^{13}. \end{eqnarray*}

\begin{figure}[h]
\hspace{0.8cm}
\epsfxsize=2.2in \epsffile{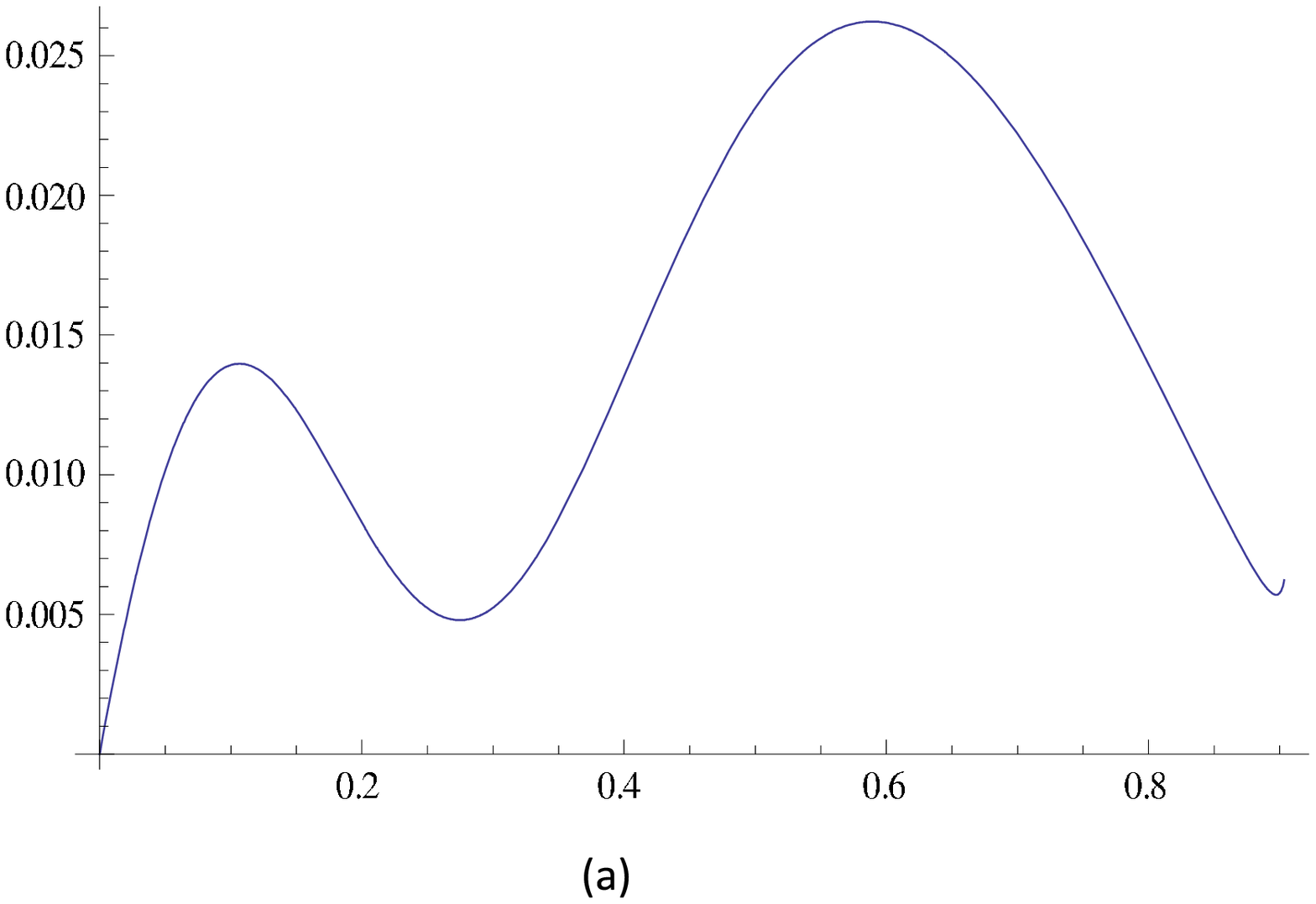}
\hspace{0.2cm}
\epsfxsize=2.2in \epsffile{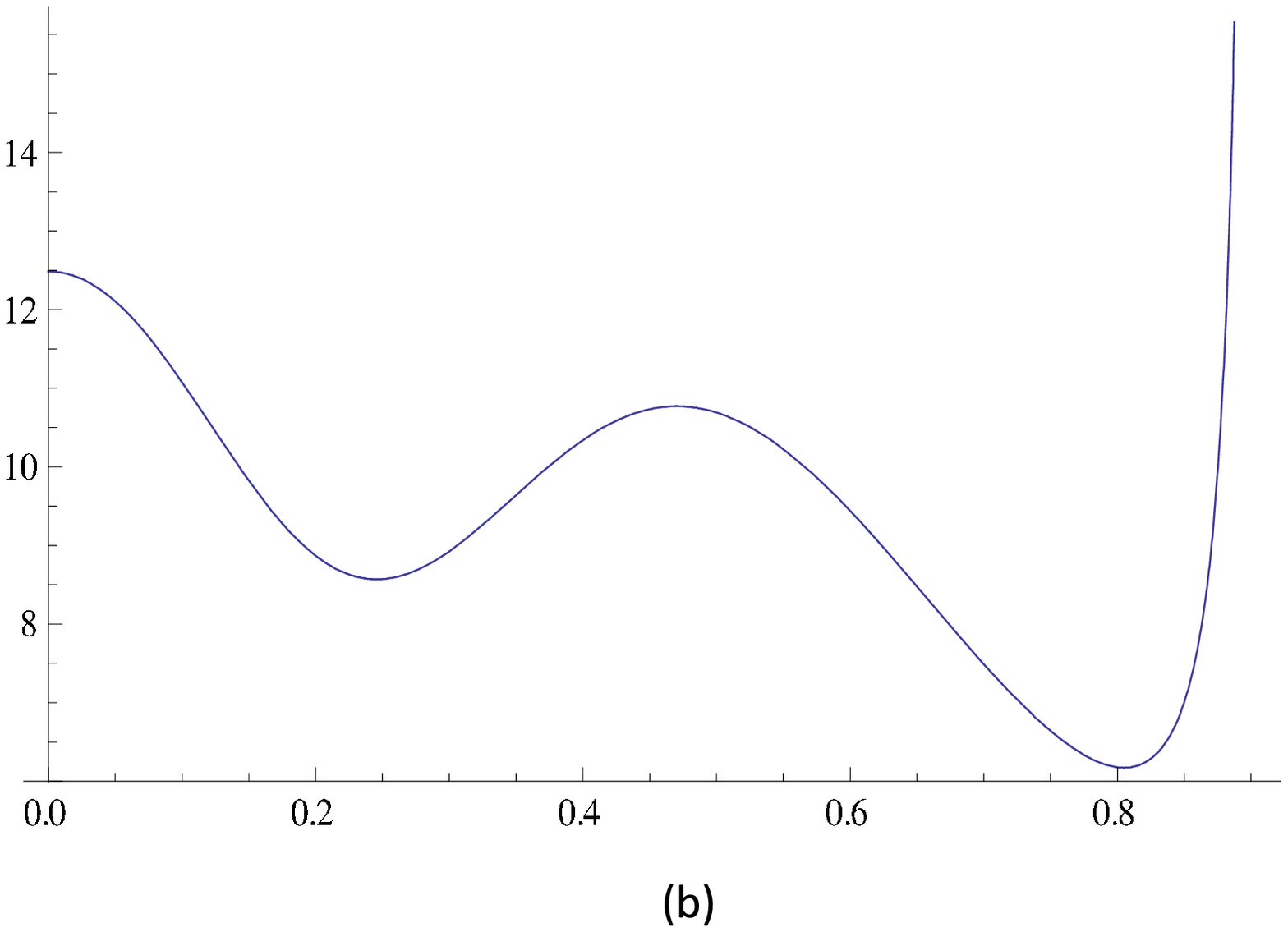}
\caption{Graphs from $\rho_X(1)-\rho_{X^2}(1)$
{\sf (a)} and $K_{p12}(r)$  {\sf (b)}, with
$0<r<\sqrt[4]{\frac{2}{3}}$\label{fig3}}
\end{figure}

As can be seen in Figure~\ref{fig3}(a),  the Taylor property is now
achieved for all considered parameterizations of Model (\ref{eq3}).

Concerning the kurtosis of this model, it is given by
$$K_{p12}(r)=\frac{-2 (-5 + 6 r^2)}{49 (-3 + 4 r^3) (-2 + 3 r^4) }\frac{N_{p12}^*(r)}{D_{p12}^*(r)}-3,$$ where
\begin{eqnarray*}
N_{p12}^*(r) & = & 599933276250 - 2617890660000 r +
      4970166270300 r^2 \\ & \; & \hspace*{0.1cm}- 5546727078200 r^3
+ 59041720498845 r^4 -
      161234870633760 r^5 \\ & \; & \hspace*{0.15cm}+ 126074334149694 r^6 + 2238307939140 r^7+
      25296348317400 r^8 \\ & \; & \hspace*{0.2cm} - 57875913071352 r^9 - 89078826937116 r^{10}+
      180941306693040 r^{11} \\  & \; & \hspace*{0.25cm}- 102607682886720 r^{12} +
      19713391884288 r^{13} \\
D_{p12}^*(r) & = &  (36300 -
     79200 r + 219255 r^2 - 171160 r^3 + 29472 r^4)^2.
\end{eqnarray*}

\vspace*{0.2cm} \noindent{\bf Error process with  Pareto  density
$\displaystyle
f(x)=\frac{9\alpha^{9}}{x^{10}}\I_{]\alpha,+\infty[}(x)$}

\noindent We have

$\beta^4 \, \mu_4<1 \;\Longleftrightarrow \;0<r
<\sqrt[4]{\frac{5}{9}}\simeq 0.863$ and
\begin{eqnarray*}
\rho_X(1) & = & \frac{8 r (15680 - 27720 r + 39564 r^2 - 27864 r^3 + 6561
r^4)}{47040 - 105840 r + 343119 r^2 - 315504 r^3 + 73791 r^4}\\
\rho_{X^2}(1)  & = &  \frac{r}{48}\,\frac{N_{p9}(r)}{D_{p9}(r)},\end{eqnarray*}
\noindent with
\begin{eqnarray*}
 N_{p9}(r) & = & -67737600 - 83339200 r + 19038600 r^2 -
88401600 r^3 \\ & \; & \hspace*{0.1cm} -
 148138920 r^4 - 511287075 r^5
+ 1466330040 r^6 + 1499354145 r^7 \\ & \; & \hspace*{0.2cm} -
 1537629480 r^8 - 1966005837 r^9 - 602608896 r^{10}
 \\
& \; & \hspace*{0.3cm} + 3869347563 r^{11} -
 61620912 r^1{2} - 2818841796 r^1{3} + 1179090432 r^{14} \\
D_{p9}(r) & = & -627200 + 235200 r - 1650600 r^2 - 8601600 r^3 - 13809280
r^4 \\ & \; & \hspace*{0.1cm} +
 31729095 r^5
+ 27010080 r^6 - 23002305 r^7 - 21773448 r^8 \\ & \; & \hspace*{0.2cm} -
 24182469 r^9 + 58517640 r^{10}  + 9248823 r^{11} \\ & \; & \hspace*{0.3cm} - 50143536 r^{12} +
 19665504 r^{13}. \end{eqnarray*}

\begin{figure}
\hspace{0.8cm}
\epsfxsize=2.2in \epsffile{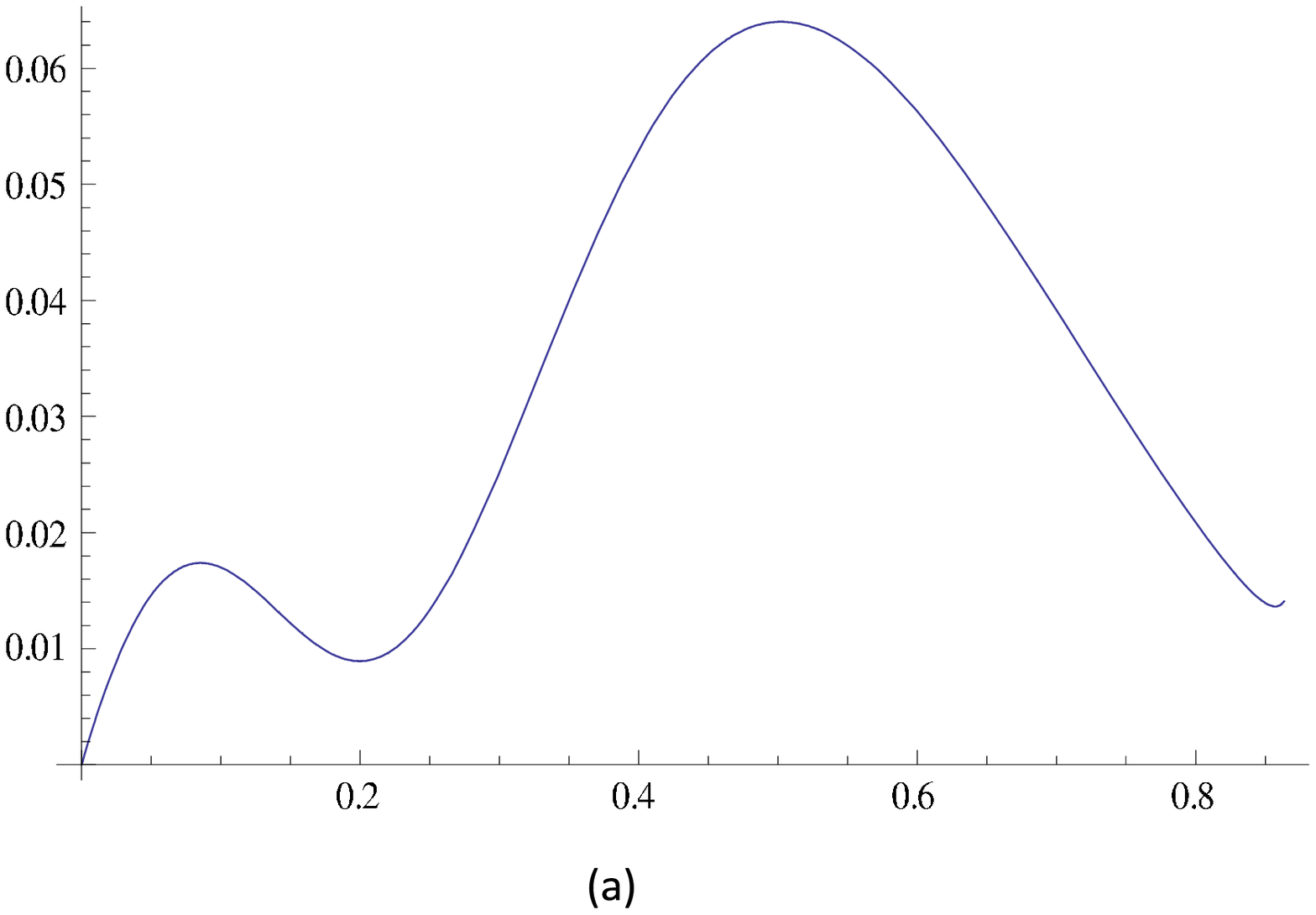}
\hspace{0.2cm}
\epsfxsize=2.2in \epsffile{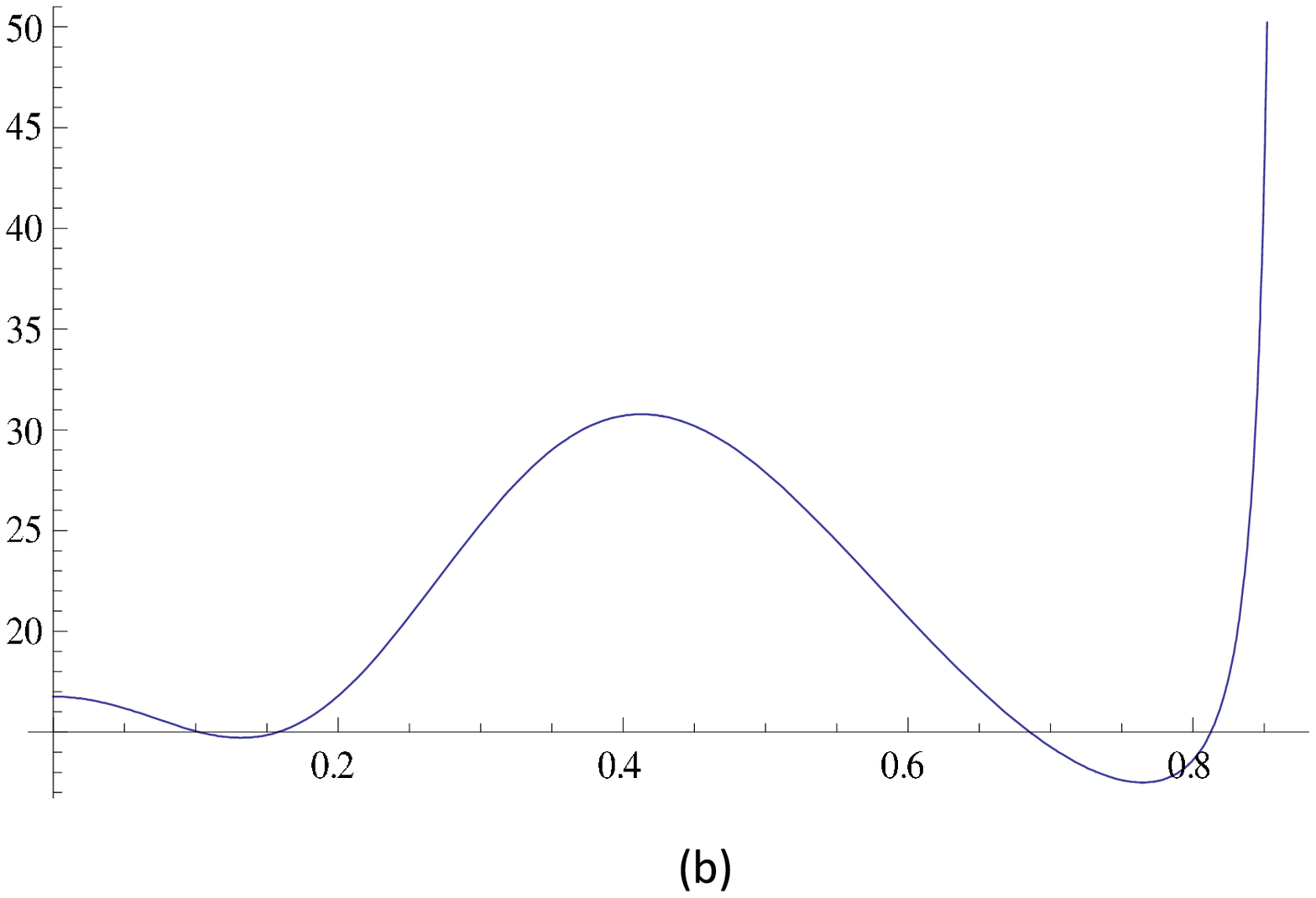}
\caption{Graphs from $\rho_X(1)-\rho_{X^2}(1)$
{\sf (a)} and $K_{p9}(r)$  {\sf (b)}, with
$0<r<\sqrt[4]{\frac{5}{9}}$\label{fig4}}
\end{figure}

 The Taylor property is also present for all considered
parameterizations of Model (\ref{eq3}), as it is illustrated in
Figure~\ref{fig4}(a), and we point out that the magnitude of the
difference $\rho_X(1)-\rho_{X^2}(1)$ is greater in this case than in
the case $\nu=12$.

The kurtosis of Model (\ref{eq3}) is now given by
$$K_{p9}(r)=\frac{7 - 9 r^2}{9 (-2 + 3 r^3) (-5 + 9 r^4)}\frac{N_{p9}^*(r)}{D_{p9}^*(r)}-3,$$ where
\begin{eqnarray*}
N_{p9}^*(r) & = & 62449049600 - 281020723200 r + 532657440000 r^2 -
      582241598400 r^3 \\
 & \; & \hspace*{0.1cm}+ 25718506014670 r^4 - 92872063045440 r^5 +
      100396353649230 r^6\\  & \; & \hspace*{0.15cm}- 6337711636725 r^7
- 8536591340550 r^8 -
      41782534519365 r^9 \\  & \; & \hspace*{0.2cm} - 62336742758694 r^{10}+
      195729014255481 r^{11}\\ & \; & \hspace*{0.25cm} - 145385404543008 r^{12} +
      35664808109193 r^{13}\\
D_{p9}^*(r) & = &  (15680 -
      35280 r + 114373 r^2 - 105168 r^3 + 24597 r^4)^2.
\end{eqnarray*}

\vspace*{5mm}
%\vspace*{5mm}

We observe that the kurtosis of the process $X$ is greater when
$\nu=9$ than when $\nu=12$, corresponding to an analogous relation
between the kurtosis of the respective error processes. In these two
examples, it is seen again how  the Taylor property emerges when the
process $X$ is leptokurtic.

\begin{figure}[h]
\centering
\epsffile{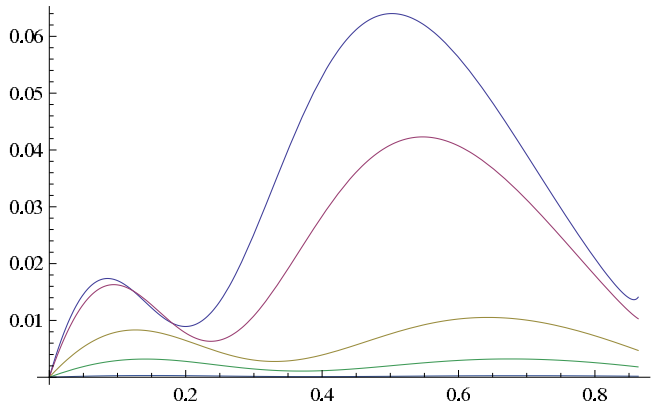}
\caption{Graphs from
$\rho_X(1)-\rho_{X^2}(1)$, $\nu=9,10,20,50,100$ (from top to
bottom), $0<r<\sqrt[4]{\frac{5}{9}}$\label{fig5}}
\end{figure}

As regards the Pareto distribution, graphic representations for
several values of $\nu$ suggest that the difference
$\rho_X(1)-\rho_{X^2}(1)$ tends to zero as $\nu$ tends to infinity
(corresponding to decreasing values of the kurtosis of the Pareto
distribution). This situation is illustrated in Figure~\ref{fig5}
and strongly contributes to conjecture that  the Taylor property and
leptokurtosis are highly related in time series.

\vspace*{5mm}
\vspace*{5mm}

\section{The Taylor property in the case of symmetrically
distributed errors: simulation study}
 \label{Sec5}

When the errors are symmetrically distributed, the autocorrelation
function of $X^2$ for model (\ref{eq1}) verifies $\rho_{X^2}(1)=0$,
if $k>1$ (Martins, [6]). So, in this case, the property
$\rho_{|X|}(1)>\rho_{X^2}(1)$ is equivalent to $\rho_{|X|}(1)>0$.
However, the autocorrelation function  of the process $(|X_t|, t \in
\Z)$ is not available when the error process is allowed to assume
negative values. To investigate the presence of  the Taylor property
in Model (\ref{eq3}) with symmetrically distributed errors, we
perform a simulation study considering the simple first-order
bilinear diagonal model with an i.i.d. error process
$(\varepsilon_t,t\in\Z)$ with four symmetrical distributions with
unit variance, namely, the uniform distribution in
$]-\sqrt{3},\sqrt{3}[$, the standard normal distribution, and the
distribution of a variable
$\varepsilon=\sqrt{\frac{\nu-2}{\nu}}\,Y$, where $Y$ has a Student
distribution with $\nu$ degrees of freedom ($\nu=30$ and $\nu=9$).
In each case, the condition $E(|\ln |\varepsilon_t||)<+\infty$ is
satisfied and parameterizations that satisfy $\beta^4\mu_4<1$ are
considered in the simulations. For each value of the parameter
$\beta$ and each one of the considered distributions, we generate
500 observations according to the corresponding model and obtain the
$95\%$ confidence intervals for the probability that such a model
satisfies  the Taylor property. The results appear in Table
\ref{tab1} (where NA means ``Not Applicable", due to the fact that
the corresponding value of $\beta$ does not satisfy the condition
$\beta^4\mu_4<1$). The special values 0.69, 0.74, 0.75 and 0.863 are
the greatest values of $\beta$ such that $\beta^4\mu_4<1$ for each
one of the considered distributions.

\vspace*{5mm}
\begin{table}[h]
\centering
\begin{tabular}{|c||c|c|c|c|} \hline
$\beta$ & $U\left(]-\sqrt{3},\sqrt{3}[\right)$ & $N(0,1)$ &
$\sqrt{\frac{14}{15}}\,Y, \;\,Y\sim T(30)$
 & $\sqrt{\frac{7}{9}}\,Y,\;\,Y\sim T(9)$\\
\hline \hline
0.01 & [0.373,0.627] & [0.459,0.708] & [0.459,0.708]  & [0.476,0.724]   \\
\hline
0.05 & [0.357,0.610] & [0.373,0.627] &  [0.373,0.627] & [0.407,0.660] \\
\hline
0.1 &  [0.140,0.360]  & [0.292,0.541] & [0.214,0.453] & [0.260,0.506] \\
\hline
0.2 &  [0,0]  & [0,0.105] & [0,0.049]  & [0,0.049]  \\
\hline
0.3 &  [0,0] &  [0,0] &  [0,0]  & [0,0.079] \\
\hline
0.4 &  [0,0] & [0,0] & [0,0.079] & [0.260,0.506] \\
\hline
0.5 & [0,0] & [0.155,0.379] &  [0.292,0.541]  & [0.699,0.901] \\
\hline
0.6 & [0,0] & [0.566,0.801] &   [0.603,0.831] & [0.781,0.953] \\
\hline
0.69 & [0,0]& [0.802,0.965] &   [0.802,0.965] & [0.951,1] \\
\hline
0.74 & [0,0.079] & [0.847,0.987]   &   [0.870,0.996] & NA \\
\hline
0.75 & [0.004,0.130] & [0.847,0.987] &   NA & NA \\
\hline
0.863&  [0.566,0.801]  & NA &   NA & NA \\
\hline
\end{tabular}\caption{$95\%$ confidence intervals for the probability that the
model with symmetrical innovations presents  the Taylor property.\label{tab1}}
\end{table}
\vspace*{5mm}

We can observe that  the Taylor property seems to be present for
high values of $\beta$ and that this presence increases with the
kurtosis of the error process, as we have established and observed
in non-negative bilinear models.

The confidence intervals corresponding to small values of $\beta$ do
not allow us to infer about the presence of  the Taylor property, as
they certainly correspond to values of $\beta$ for which the
difference $\rho_X(1)-\rho_{X^2}(1)$ is close to zero.

\section{Conclusions}
 \label{Sec5}

The studies presented here show that bilinear models are able to
reproduce  the Taylor effect. They also reinforce the connection of
the Taylor property to leptokurtic models which has been observed in
the few theoretical studies developed until now. In fact, He and
Ter\"asvirta (\cite{HeTer}), Gon\c{c}alves, Leite and Mendes-Lopes
(\cite{GLML}) and Haas (\cite{Haas}) show the presence of this
property in some conditional heteroskedastic models, which are
leptokurtic processes. Moreover, all the cases considered in this
paper, also show that, when  the Taylor property occurs, the model
is leptokurtic.

We still observe that leptokurtosis is not enough to induce the
Taylor property. Examples of bilinear models that are leptokurtic
but do not have  the Taylor property are
$X_t=X_{t-1}\varepsilon_{t-1}+\varepsilon_t$, where $\varepsilon_t$
is uniformly distributed in $[0,1]$, and $X_t=0.5
X_{t-1}\varepsilon_{t-1}+\varepsilon_t$, where $\varepsilon_t$ is
exponentially distributed with mean $0.2$. This is in line with the
simulation results of He and Ter\"asvirta (\cite{HeTer}) suggesting
that  the Taylor property is not present for the standard
GARCH$(1,1)$ process with normal errors.

In conclusion, our study allows to conjecture that a general
assessment of  the Taylor property in the bilinear process is
strongly dependent on its tails weight.

\end{document}